\documentclass[a4paper]{article}
\usepackage{amssymb,amsmath,amsthm}
\usepackage{graphicx}
\usepackage{mathrsfs}
\usepackage{hyperref}

\theoremstyle{definition}
\newtheorem{thm}{Theorem}[section]

\newtheorem{lem}{Lemma}[section]
\theoremstyle{remark}

\numberwithin{equation}{section}

\newcommand{\scurl}{\operatorname{curl}}
\newcommand{\curl}{\operatorname{\text{\bf curl}}}
\newcommand{\bs}{\boldsymbol}
\newcommand{\wt}{\widetilde}
\newcommand{\eps}{^\epsilon}
\newcommand{\loc}{\text{loc}}
\newcommand{\re}{\text{Re}\,}
\newcommand{\pd}[2]{\dfrac{\partial #1}{\partial #2}}
\newcommand{\build}[3]{\mathrel{\mathop{\kern 0pt#1}\limits_{#2}^{#3}}}

\def\diver{\operatorname{div}}

\begin{document}

\title{A Two-dimensional eddy current model using thin inductors}
\author{Youcef {\sc Amirat}, Rachid {\sc Touzani}\\
\small Laboratoire de Math\'ematiques, UMR CNRS 6620\\
\small Universit\'e Blaise Pascal (Clermont--Ferrand)\\
\small 63177 Aubi\`ere cedex, France\\
\small \tt \{Youcef.Amirat,Rachid.Touzani\}@math.univ-bpclermont.fr}
\date{}

\maketitle

\begin{abstract}
We derive a mathematical model for eddy currents in two dimensional geometries where
the conductors are thin domains. We assume that the current flows in the $x_3$-direction
and the inductors are domains with small diameters of order $O(\epsilon)$.
The model is derived by taking the limit $\epsilon\to 0$.
A convergence rate of $O(\epsilon^\alpha)$ with $0<\alpha<1/2$ in the $L^2$--norm is shown
as well as weak convergence in the $W^{1,p}$ spaces for $1<p<2$.

\end{abstract}


\section{Introduction}
Mathematical modelling of eddy current problems often involves multiple conductors
with various sizes. Typically, electrotechnical devices involve thin conductors as wires
or coils as well as massive conductors. Numerical solution of such problems may then
encounter serious difficulties in the choice of the domain meshes which can in particular
lead to ill conditioning. Asymptotic analysis of these problems appears as an efficient tool
to obtain limit problems that are simpler to solve and better conditioned.

We consider, in the present work, a two-dimensional eddy current problem, formulated in
terms of a scalar potential in the whole plane. The electrically conducting domain consists
in a ``thick" conductor $\Omega$ and two ``thin" domains assumed to carry the same current
with opposite sign. In terms of the current conservation principle, this means that these
inductors are assumed to be virtually \emph{linked at the infinity}. The derivation of the
model for thin inductors is obtained by assuming that these domains are of small diameters of
order $\epsilon\ll 1$. We show that taking the limit when these
diameters tend to zero leads to a singular elliptic problem, the singularity being due to the presence
of Dirac measures.

The outline of the paper is the following: We start in Section 2 by deriving the considered
eddy current model from a 3-D model. We emphasize on a careful modelling that takes into
account the total current flowing in the inductors. In Section 3, we state the main convergence
result and prove it through some preliminary lemmas. Section 4 is devoted to further convergence
results in $W^{1,p}$ spaces.

\section{Statement of the problem}

Let $\Lambda=\Omega\times\mathbb R$ denote a cylindrical conductor where $\Omega$ is a domain in $\mathbb R^2$
with a smooth boundary $\Gamma$. We assume that the domain $\Omega$ is the union of three connected domains
$\Omega_k$ with respective boundaries $\Gamma_k$, $k=0,1,2$ (see Figure 1), and that the
closures of the domains $\Omega_k$ are disjointed.
We shall also deal with the complement $\Omega'=\mathbb R^2\setminus\overline\Omega$ of $\Omega$.

\begin{figure}[!ht]
\includegraphics[bb=100 580 596 680]{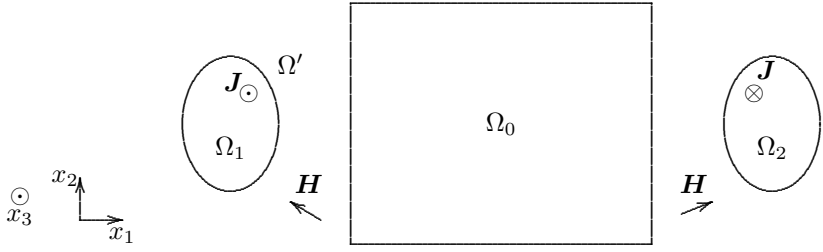}%
\label{fig}\caption{A typical configuration of the conductors}
\end{figure}

In the following we shall make use of generic constants that do not depend on the small parameter $\epsilon$.
Time harmonic eddy currents equations read:
\begin{equation}
\begin{aligned}{}
&\curl\bs H - \bs J = 0&&\qquad\text{in }\Lambda,\\
&\curl\bs H = 0&&\qquad\text{in }\mathbb R^3\setminus\overline\Lambda,\\
&i\omega\mu\bs H + \curl\,(\sigma^{-1}\bs J) = 0&&\qquad\text{in }\Lambda,\\
&\diver\,(\mu\bs H) = 0&&\qquad\text{in }\mathbb R^3.
\end{aligned}
\label{EC-3D}
\end{equation}
Here the vector fields $\bs H$ and $\bs J$ denote respectively the magnetic field and
the current density. Moreover, $\sigma$ and $\mu$ are respectively the electric conductivity and the magnetic permeability.
We assume for the sake of simplicity that $\sigma$ and $\mu$ are positive constants.
In order to take advantage of the geometry of $\Lambda$, we seek unknowns in the form:
\begin{align}
&\bs H(x_1,x_2,x_3) = H_1(x_1,x_2)\,\bs e_1 + H_2(x_1,x_2),0)\,\bs e_2,\nonumber\\
&\bs J(x_1,x_2,x_3) = J(x_1,x_2)\,\bs e_3,\nonumber
\end{align}
where $(\bs e_1,\bs e_2,\bs e_3)$ is the canonical basis of $\mathbb R^3$.
Equations \eqref{EC-3D} become then
\begin{alignat}{2}
&\scurl\bs H - J = 0                          &&\qquad\text{in }\Omega,            \label{EC-1}\\
&\scurl\bs H = 0                              &&\qquad\text{in }\Omega',           \label{EC-2}\\
&i\omega\mu\bs H + \curl\,(\sigma^{-1}J) = 0  &&\qquad\text{in }\Omega,            \label{EC-3}\\
&\diver\,\bs H = 0                            &&\qquad\text{in }\mathbb R^2.       \label{EC-4}
\end{alignat}
where $\scurl$ and $\curl$ denote respectively the scalar and vector curl operator in
2-D, \emph{i.e.}
$$
\scurl\bs u := \pd{H_2}{x_1} - \pd{H_1}{x_2},\quad
\curl\varphi = \pd{\varphi}{x_2}\,\bs e_1 - \pd{\varphi}{x_1}\,\bs e_2.
$$
Our aim now is to derive a simple model for eddy currents. Using Equation \eqref{EC-4},
we deduce the existence of a scalar potential $u:\mathbb R^2\to\mathbb C$ such that
\begin{equation}
\mu\bs H = \curl u\qquad\text{in }\mathbb R^2.\label{divB=0}
\end{equation}
Equation \eqref{EC-1} yields
$$
\begin{aligned}{}
&\scurl\,\curl u = \mu J&&\qquad\text{in }\Omega,\\
&\scurl\,\curl u = 0&&\qquad\text{in }\Omega',
\end{aligned}
$$
or equivalently,
\begin{equation}
\begin{aligned}{}
-&\Delta u = \mu J&&\qquad\text{in }\Omega,\\
&\Delta u = 0&&\qquad\text{in }\Omega'.\label{eq1}
\end{aligned}
\end{equation}
On the other hand we obtain from \eqref{EC-3} and \eqref{divB=0},
$$
\curl\,(i\omega u + \sigma^{-1}J) = 0.
$$
Whence
\begin{equation}
i\omega\sigma u + J = \sigma C_k\qquad\text{in }\Omega_k,\label{eq-Ck}
\end{equation}
where $C_k$ are complex constants, for $k=0,1,2$. Replacing this in \eqref{eq1}, we obtain
$$
\begin{aligned}{}
&-\Delta u + i\omega\mu\sigma u = \mu\sigma C_k&&\qquad\text{in }\Omega_k,\ k=0,1,2,\\
&\Delta u = 0&&\qquad\text{in }\Omega'.
\end{aligned}
$$
Finally, various considerations dealing with interface conditions and the behaviour at the infinity lead to the
problem:
\begin{equation}
\begin{aligned}{}
&-\Delta u + i\omega\mu\sigma u = \mu\sigma C_k  &&\qquad\text{in }\Omega_k,\ k=0,1,2,\\
&\Delta u = 0                                    &&\qquad\text{in }\Omega',\\
&[u] = 0                                         &&\qquad\text{on }\Gamma,\\
&\Big[\pd un\Big] = 0                            &&\qquad\text{on }\Gamma,\\
&u(x) = \alpha + O(|x|^{-1})                     &&\qquad\text{as } |x|\to\infty.
\end{aligned}
\end{equation}
Here above, $[\,\cdot\,]$ denotes the jump of a function across the boundary $\Gamma$, this jump being equal to
the external trace minus the internal one.
It remains to determine the constants $C_k$ in function of problem data. For this end,
it turns out to be realistic to prescribe the total current in each conductor. We then
assume that this quantity, denoted by $I$ is given as
\begin{equation}
\label{I}
\int_{\Omega_1} J\,dx = -\int_{\Omega_2} J\,dx = I.
\end{equation}
Note that the first identity is imposed in order to enforce a current conservation principle.
For the same reason, we impose
$$
\int_{\Omega_0} J\,dx = 0.
$$
Making use of these conditions, we obtain for the constants $C_k$, the values
\begin{align*}
&C_0 = i\omega\wt u_1,\\
&C_1 = i\omega\wt u_1 + \frac I{\sigma\,|\Omega_1|},\\
&C_2 = i\omega\wt u_2 - \frac I{\sigma\,|\Omega_2|},
\end{align*}
where $|\Omega_k|$ stands for the measure of $\Omega_k$ and $\wt u_k$ is the average of $u$
on $\Omega_k$, \emph{i.e.}
$$
\wt u_k := \frac 1{|\Omega_k|}\int_{\Omega_k}u\,dx.
$$
We obtain the problem:
\begin{equation}
\begin{aligned}{}
&-\Delta u + i\omega\mu\sigma(u-\wt u_0) =  0                          &&\qquad\text{in }\Omega_0,\\
&-\Delta u + i\omega\mu\sigma(u-\wt u_1) =  \frac {\mu I}{|\Omega_1|}  &&\qquad\text{in }\Omega_1,\\
&-\Delta u + i\omega\mu\sigma(u-\wt u_2) = -\frac {\mu I}{|\Omega_2|}  &&\qquad\text{in }\Omega_2,\\
&\Delta u = 0                                                          &&\qquad\text{in }\Omega',\\
&[u] = 0                                                               &&\qquad\text{on }\Gamma,\\
&\Big[\pd un\Big] = 0                                                  &&\qquad\text{on }\Gamma,\\
&u(x) = \alpha + O(|x|^{-1})                                           &&\qquad\text{as } |x|\to\infty.
\end{aligned}\label{Pb-I}
\end{equation}
Note that, owing to \eqref{I}, the solution of Problem \eqref{Pb-I} is known up to an additive constant. For this reason,
we impose the condition
\begin{equation}
\wt u_0=0,\label{ZeroAverage}
\end{equation}
which enforces a value for the constant $\alpha$.

Let us prove that Problem \eqref{Pb-I}--\eqref{ZeroAverage} has a unique solution.
We define, for this end, the Beppo-Levi space (see \cite{Nedelec}),
$$
W^1(\mathbb R^2) := \Big\{v;\ \rho\,v\in L^2(\mathbb R^2),\ \nabla v
\in L^2(\mathbb R^2)^2\Big\},
$$
where $\rho$ is the weight function given by
\begin{equation}
\rho(x)=\frac 1{(1+|x|)\log(2+|x|)}.\label{rho}
\end{equation}
We furthermore define the space
$$
V := \{v\in W^1(\mathbb R^2);\ \wt v_0=0\}.
$$
It is well known (cf. \cite{RT}) that the semi-norm
$$
|v|_{W^1(\mathbb R^2)} := \Big(\int_{\mathbb R^2}|\nabla v|^2\,dx\Big)^{\frac 12}
$$
is a norm on the space $V$, equivalent to the one induced by $W^1(\mathbb R^2)$, \emph{i.e.}
we have in particular
\begin{equation}
\|\rho\, v\|_{L^2(\mathbb R^2)} \le C\,|v|_{W^1(\mathbb R^2)}\qquad\forall\ v\in V.\label{EqNorm}
\end{equation}
Here and in the following $|\nabla v|$ stands for the function
$$
|\nabla v| = \Big(\pd v{x_1}\pd {\overline v}{x_1}+\pd v{x_2}\pd {\overline v}{x_2}\Big)^{\frac 12}.
$$

A variational formulation of \eqref{Pb-I} consists in seeking a function $u\in V$ such that
\begin{equation}
\begin{aligned}{}
\int_{\mathbb R^2}\nabla u\cdot\nabla\overline v\,dx + i\beta\sum_{k=0}^2\int_{\Omega_k}(u-\wt u_k)\overline v\,dx
= \mu I(\overline{\wt v}_1-\overline{\wt v}_2)\quad\forall\ v\in V,
\end{aligned}
\label{pb-I-var}
\end{equation}
where $\beta=\omega\mu\sigma$ and $\overline v$ is the complex conjugate of $v$.

\begin{thm}
Problem \eqref{pb-I-var} has a unique solution.
\end{thm}

\begin{proof}
Let us define, for $u,v\in V$, the sesquilinear and antilinear forms,
\begin{align*}
&a(u,v) := \int_{\mathbb R^2}\nabla u\cdot\nabla\overline v\,dx +
i\beta\sum_{k=0}^2\int_{\Omega_k}(u-\wt u_k)\overline v\,dx,\\
&L(v) := \mu I(\overline{\wt v}_1-\overline{\wt v}_2).
\end{align*}
The forms $a$ and $L$ are obviously continuous. In addition, since
$$
\int_{\Omega_k}(u-\wt u_k)\overline v\,dx = \int_{\Omega_k}(u-\wt u_k)(\overline v-\overline{\wt v}_k)\,dx,
$$
then we have
\begin{align*}
a(v,v) &= \int_{\mathbb R^2}|\nabla v|^2\,dx +
i\beta\sum_{k=0}^2\int_{\Omega_k}(v-\wt v_k)\overline v\,dx\\
&= \int_{\mathbb R^2}|\nabla v|^2\,dx +
i\beta\sum_{k=0}^2\int_{\Omega_k}|v-\wt v_k|^2\,dx.
\end{align*}
Then
$$
\re(a(v,v)) = \int_{\mathbb R^2}|\nabla v|^2\,dx.
$$
We deduce then that $a$ is coercive on $V$ and the Lax-Milgram theorem gives the existence
and uniqueness of a solution $u\in V$ to \eqref{pb-I-var}.
\end{proof}

We now consider that the domains $\Omega_1$ and $\Omega_2$ are thin in the following sense:
we define the domain $\Omega\eps_k:=\Omega_k$ by
$$
\Omega\eps_k = z_k + \epsilon\,\widehat\Omega_k\qquad k=1,2,
$$
where $\epsilon$ is a small positive number, $z_k\in\mathbb R^2$, and $\widehat\Omega_k$ is a smooth domain in $\mathbb R^2$.
We assume furthermore that the domains $\overline\Omega_k\eps$ and $\overline\Omega_0$
are disjointed for $\epsilon$ small enough. Furthermore, we denote in the following by $\Omega\eps$
the union $\Omega_0\cup\Omega_1\eps\cup\Omega_2\eps$. Finally, let us mention that, throughout this
paper, $C, C_1, C_2, \ldots$ will stand for generic constants that do not depend on $\epsilon$.
Our aim is to study the asymptotic behavior, as $\epsilon\to 0$, of the solution $u$ to Problem \eqref{Pb-I}.

\section{The limit problem}

Let us first, for clarity, rewrite Problem \eqref{Pb-I} with the parameter $\epsilon$.
Denoting by $\chi_0$ and $\chi_k\eps$ the characteristic functions of $\Omega_0$
and $\Omega_k\eps$, respectively, we have
\begin{equation}
\begin{aligned}{}
&-\Delta u\eps + i\beta\sum_{k=1}^2\chi_k\eps\,(u\eps-\wt u\eps_k)
+ i\beta\chi_0\, u\eps = \mu I\Big(\frac{\chi_1\eps}{|\Omega_1\eps|}
-\frac{\chi_2\eps}{|\Omega_2\eps|}\Big)&&\text{in }\mathbb R^2,\\
&u\eps(x) = \alpha + O(|x|^{-1})&&\text{as } |x|\to\infty.
\end{aligned}\label{Pb-Ie}
\end{equation}
Let us recall that the condition \eqref{ZeroAverage} fixes the value of $\alpha$.

We next define the weighted space that will be used for convergence results:
$$
L^2_\rho(\mathbb R^2) = \{v;\ \rho v\in L^2(\mathbb R^2)\}.
$$

We also define a problem that will be defined as the limit problem. This one is the following:
\begin{equation}
\begin{aligned}{}
&-\Delta u + i\beta\chi_0\,u = \mu I(\delta_{z_1}-\delta_{z_2})&&\ \text{in }\mathbb R^2,\\
&u(x) = \alpha + O(|x|^{-1})&&\text{ as } |x|\to\infty,
\end{aligned}\label{Pb-Lim}
\end{equation}
where $\delta_{z_k}$ is the Dirac measure concentrated at $z_k$.
For Problem \eqref{Pb-Lim}, we need a uniqueness result. Let us define for this the notion
of weak solution. We shall say in the sequel that $u$ is a weak $L^2_\rho$--solution
of Problem \eqref{Pb-Lim} if $u\in L^2_\rho(\mathbb R^2)$ and if we have
\begin{equation}
\int_{\mathbb R^2} u\,(-\Delta\overline\varphi + i\beta\chi_0\,\overline\varphi)\,dx =
\mu I\,(\overline\varphi(z_1)-\overline\varphi(z_2))
\qquad\forall\ \varphi\in \mathscr D(\mathbb R^2),\label{WeakSol}
\end{equation}
where $\mathscr D(\mathbb R^2)$ is the space of indefinitely differentiable functions
with compact support in $\mathbb R^2$.

\begin{lem}
\label{Uniqueness}
Problem \eqref{Pb-Lim} has at most one weak $L^2_\rho$--solution.
\end{lem}

\begin{proof}
Let $u_1$ and $u_2$ denote two weak $L^2_\rho$--solutions of
\eqref{WeakSol}. The difference $u=u_1-u_2$ satisfies then
$$
\int_{\mathbb R^2} u\,(-\Delta\overline\varphi + i\beta\chi_0\,\overline\varphi)\,dx = 0
\qquad\forall\ \varphi\in\mathscr D(\mathbb R^2).
$$
This relation is still true for all functions $\varphi\in L^2_\rho(\mathbb R^2)$ with
$$
\int_{\mathbb R^2}\rho^2 u\,\overline\psi\,dx = 0\qquad\forall\ \psi\in L^2_\rho(\mathbb R^2),
$$
where
\begin{equation}
-\Delta\overline\varphi + i\beta\chi_0\,\overline\varphi = \rho^2\overline\psi\qquad\text{in }\mathbb R^2.
\label{weak-dual}
\end{equation}
Note that Equation \eqref{weak-dual} admits a unique solution in $W^1(\mathbb R^2)$.
Choosing $\psi=u$, we deduce
$$
\int_{\mathbb R^2}\rho^2\,|u|^2\,dx = 0.
$$
This implies $u=0$ and uniqueness follows.
\end{proof}

We now state the first convergence result.

%
%

\begin{thm}
\label{CV}
The sequence $(u\eps)$ converges in $L^2_\rho(\mathbb R^2)$, when $\epsilon\to 0$, to the unique solution of Problem
\eqref{Pb-Lim}.
\end{thm}

The remaining of this section is devoted to the proof of Theorem \ref{CV}.
It is clear that the structure of the right-hand side in Problem \eqref{Pb-Ie} suggests
that the convergence cannot be obtained in the space $W^1(\mathbb R^2)$.
To obtain a weaker result we resort to a duality technique due to Lions-Magenes (\cite{LM}, p. 177)
and Damlamian-Ta Tsien Li \cite{DL}.

Let, in the following, $B$ denote a ball that contains the domains $\overline\Omega_0$,
$\overline\Omega_1\eps$ and $\overline\Omega_2\eps$ for all $\epsilon\ll 1$.
Multiplying Equation \eqref{Pb-Ie} by a test function $\varphi\in V\cap H^2_\loc(\mathbb R^2)$,
and using the Green formula, we obtain
$$
-\int_{\mathbb R^2}u\eps\Delta\overline\varphi\,dx +
i\beta\sum_{k=1}^2\int_{\Omega_k\eps}(u\eps-\wt u\eps_k)\overline\varphi\,dx
+ i\beta\int_{\Omega_0}u\eps\overline\varphi\,dx =
\mu I(\overline{\wt\varphi}_1-\overline{\wt\varphi}_2).
$$
Since
\begin{equation}
\int_{\Omega_k\eps}(u\eps-\wt u\eps_k)\,\overline\varphi\,dx =
\int_{\Omega_k\eps}u\eps(\overline\varphi-\overline{\wt\varphi}_k\eps)\,dx
=\int_{\Omega_k\eps}(u\eps-\wt u\eps_k)\,(\overline\varphi-\overline{\wt\varphi}_k\eps)\,dx,
\label{Ident}
\end{equation}
we deduce that
\begin{equation}
\int_{\mathbb R^2}u\eps\Big(-\Delta\overline\varphi + i\beta\sum_{k=1}^2\chi_k\eps(\overline\varphi
-\overline{\wt\varphi}\eps_k) + i\beta\chi_0\overline\varphi\Big)\,dx
= \mu I(\overline{\wt\varphi}_1-\overline{\wt\varphi}_2).
\label{Adjoint-A}
\end{equation}
Let $\psi$ denote a function in $L^2_\rho(\mathbb R^2)$.
Identity \eqref{Adjoint-A} can also be written as
\begin{equation}
\int_{\mathbb R^2}\rho^2 u\eps\overline\psi\,dx = \mu I(\overline{\wt\varphi}_1\eps-\overline{\wt\varphi}_2\eps),
\label{Dual}
\end{equation}
where $\varphi\eps$ is the solution in $V\cap H^2_{\text{loc}}(\mathbb R^2)$ of
\begin{equation}
-\Delta\varphi\eps + i\beta\sum_{k=1}^2\chi_k\eps(\varphi\eps-\wt\varphi_k\eps)
+ i\beta\chi_0\varphi\eps = \rho^2\psi\qquad\text{in }\mathbb R^2.
\label{Adjoint}
\end{equation}

\begin{lem}
\label{lem1}
We have the estimates:
\begin{align}
&\|\nabla\varphi\eps\|_{L^2(\mathbb R^2)^2} + \|\varphi\eps\|_{L^2(\Omega_0)} +
\epsilon^{-1}\sum_{k=1}^2\|\varphi\eps-\wt\varphi\eps_k\|_{L^2(\Omega_k\eps)} \le C\,\|\rho\psi\|_{L^2(\mathbb R^2)}.\label{Estim-phi}\\
&\|\varphi\eps\|_{H^2(B)}\le C\,\|\rho\psi\|_{L^2(\mathbb R^2)},\label{Estim-phi-a}
\end{align}
for each ball $B$ of $\mathbb R^2$ containing $\Omega\eps$.
\end{lem}

\begin{proof}
By the Green's formula, we have from \eqref{Adjoint} and Identity \eqref{Ident}
$$
\int_{\mathbb R^2}|\nabla\varphi\eps|^2\,dx + i\beta\sum_{k=1}^2\int_{\Omega_k\eps}
|\varphi\eps-\wt\varphi_k\eps|^2\,dx + i\beta\int_{\Omega_0}|\varphi\eps|^2\,dx
= \int_{\mathbb R^2}\rho^2\psi\,\overline\varphi\eps\,dx.
$$
From this and \eqref{EqNorm} we deduce that
$$
\int_{\mathbb R^2}|\nabla\varphi\eps|^2\,dx\le\|\rho\psi\|_{L^2(\mathbb R^2)}\,
\|\rho\varphi\eps\|_{L^2(\mathbb R^2)}\le C\,
\|\rho\psi\|_{L^2(\mathbb R^2)}\,\|\nabla\varphi\eps\|_{L^2(\mathbb R^2)^2},
$$
and then
\begin{equation}
\label{Estim-phi-1}
\Big(\int_{\mathbb R^2}|\nabla\varphi^\epsilon|^2\,dx\Big)^{\frac 12}\le C\,\|\rho\psi\|_{L^2(\mathbb R^2)}.
\end{equation}
Therefore, the sequence $(\varphi\eps)$ is bounded in $W^1(\mathbb R^2)$.
The $L^2$-error estimate is obtained by using the Poincar\'e-Wirtinger inequality (see \cite{B}, p.~194).
We have indeed by using \eqref{Estim-phi-1}, and since the diameter of $\Omega\eps_k$ is an $O(\epsilon)$,
\begin{align*}
&\|\varphi\eps\|_{L^2(\Omega_0)} \le C_1\,\|\nabla\varphi\eps\|_{L^2(\mathbb R^2)^2}\le C_2\,\|\rho\psi\|_{L^2(\mathbb R^2)},\\
&\|\varphi\eps-\widetilde\varphi\eps_k\|_{L^2(\Omega_k\eps)}
\le C_3\epsilon\,\|\nabla\varphi\eps\|_{L^2(\mathbb R^2)^2}
\le C_4\epsilon\,\|\rho\psi\|_{L^2(\mathbb R^2)},\qquad k=1,2.
\end{align*}
In order to prove the $H^2$--estimate, we use standard regularity results for elliptic equations
(See \cite{GT}, p.~183 for instance). We obtain for any ball $B$ of $\mathbb R^2$ containing
$\Omega\eps$, and any regular domain $D$ containing $\overline B$,
\begin{align*}
\|\varphi\eps\|_{H^2(B)}
&\le C_1\,\Big(\|\varphi\eps\|_{H^1(D)} + \|\rho^2\psi\|_{L^2(D)} +
\|\varphi\eps\|_{L^2(\Omega_0)}+
\sum_{k=1}^2\|\varphi\eps-\widetilde\varphi_k\eps\|_{L^2(\Omega_k\eps)}\Big)\\
&\le C_2\,\|\rho\psi\|_{L^2(\mathbb R^2)}.
\end{align*}
Note that the constant $C_2$ depends on the domain $B$ but does not depend on $\epsilon$.
\end{proof}

The estimates obtained in Lemma \ref{lem1} enable concluding that a subsequence
of $(\varphi\eps)$ converges toward $\varphi$ weakly in $H^2(B)$ for any ball $B$ of $\mathbb R^2$. We now characterize
the limit function.

\begin{lem}
The sequence $(\varphi\eps)$ converges, when $\epsilon\to 0$, in $W^1(\mathbb R^2)$ to
the unique solution of the equation:
\begin{equation}
-\Delta\varphi + i\beta\chi_0\varphi = \rho^2\psi\qquad\text{in }\mathbb R^2,
\label{AdjointL}
\end{equation}
Moreover, we have the error estimates
\begin{equation}
\sum_{k=1}^2\|\varphi\eps-\varphi\|_{L^2(\Omega\eps_k)} + \|\varphi\eps-\varphi\|_{H^2(B)} +
\|\nabla(\varphi\eps-\varphi)\|_{L^2(\mathbb R^2)^2} \le C\epsilon\,\|\rho\psi\|_{L^2(\mathbb R^2)},
\label{ErrorPhi}
\end{equation}
for any ball $B$ of $\mathbb R^2$ containing $\Omega\eps$.
\end{lem}

\begin{proof}
Let $\phi\eps=\varphi\eps-\varphi$. Then $\phi\eps\in V\cap H^2_\loc(\mathbb R^2)$ and satisfies the variational equation
$$
\int_{\mathbb R^2}\nabla\phi\eps\cdot\nabla\overline v\,dx + i\beta\int_{\Omega_0}\phi\eps\overline v\,dx
+ i\beta\sum_{k=1}^2\int_{\Omega_k\eps}(\varphi\eps-\widetilde\varphi_k\eps)\overline v\,dx = 0
\qquad\forall\ v\in V.
$$
Choosing $v=\phi\eps$, we obtain
$$
\int_{\mathbb R^2}|\nabla\phi\eps|^2\,dx + i\beta\int_{\Omega_0}|\phi\eps|^2\,dx
= -i\beta\sum_{k=1}^2\int_{\Omega_k\eps}(\varphi\eps-\wt\varphi_k\eps)\overline\phi\eps\,dx.
$$
Then using the estimates \eqref{Estim-phi}, we have
\begin{align*}
\int_{\mathbb R^2}|\nabla\phi\eps|^2\,dx + \int_{\Omega_0}|\phi\eps|^2\,dx
&\le \beta\,\sum_{k=1}^2\Big(\int_{\Omega_k\eps}|\varphi\eps-\wt\varphi\eps_k|^2\,dx\Big)^{\frac 12}
\Big(\int_{\Omega_k\eps}|\phi\eps|^2\,dx\Big)^{\frac 12}\\
&\le C_1\epsilon\,\|\rho\psi\|_{L^2(\mathbb R^2)}\,\|\nabla\phi\eps\|_{L^2(\mathbb R^2)^2}.
\end{align*}
Therefore, we have the bounds
\begin{align}
&\|\nabla\phi\eps\|_{L^2(\mathbb R^2)^2} \le C_1\epsilon\,\|\rho\psi\|_{L^2(\mathbb R^2)},\label{Bound-0}\\
&\|\phi\eps\|_{L^2(\Omega_0)} \le C_2\,\|\nabla\phi\eps\|_{L^2(\mathbb R^2)^2}
\le C_3\,\epsilon\,\|\rho\psi\|_{L^2(\mathbb R^2)}.\label{Bound-1}
\end{align}
The sequence $(\varphi\eps)$ converges then to $\varphi$ strongly in $W^1(\mathbb R^2)$, which yields the
limit problem \eqref{AdjointL}.

To prove the $L^2$-error estimate, we have from \eqref{EqNorm} and \eqref{Bound-0}, for $k=1,2$,
$$
\|\phi\eps\|_{L^2(\Omega_k\eps)}
\le C_1\,\|\rho\phi\eps\|_{L^2(\Omega_k\eps)}
\le C_1\,\|\rho\phi\eps\|_{L^2(\mathbb R^2)}
\le C_2\,\|\nabla\phi\eps\|_{L^2(\mathbb R^2)}
\le C_3\epsilon\,\|\rho\psi\|_{L^2(\mathbb R^2)}.
$$
The $H^2$--estimate is handled in the following way: By subtracting \eqref{Adjoint} from
\eqref{AdjointL}, we obtain
$$
-\Delta\phi\eps = -i\beta\chi_0\phi\eps - i\beta\sum_{k=1}^2\chi_k\eps(\varphi\eps-\wt\varphi_k\eps)
\qquad\text{in }\mathbb R^2.
$$
Using \eqref{Estim-phi}, \eqref{Bound-1} and classical regularity results for elliptic problems
(See \cite{GT}, p.~183 for instance), we get
\begin{align*}
\|\phi\eps\|_{H^2(B)}
&\le C_1\,\Big(\|\phi\eps\|_{H^1(D)} + \|\phi\eps\|_{L^2(\Omega_0)}+\sum_{k=1}^2\|\varphi\eps-\widetilde\varphi_k\eps\|_{L^2(\Omega_k\eps)}\Big)\\
&\le C_2\,\Big(\epsilon\,\|\nabla\phi\eps\|_{L^2(\mathbb R^2)^2} + \|\phi\eps\|_{L^2(\Omega_0)}+\epsilon\,\|\rho\psi\|_{L^2(\mathbb R^2)^2}\Big)\\
&\le C_3\,\epsilon\,\|\rho\psi\|_{L^2(\mathbb R^2)},
\end{align*}
for all compact subsets $B$ of $\mathbb R^2$ and all regular domains $D$ that contain $\overline B$.
Note that the constant $C$ depends actually on $B$.
\end{proof}

We are now ready to obtain the first convergence result for $u\eps$.
\begin{thm}
\label{ErrEst-L2}
There exists a constant $C$, independent of $\epsilon$, such that
$$
\|\rho(u-u\eps)\|_{L^2(\mathbb R^2)}\le C\epsilon^{\alpha/2}\qquad 0<\alpha<1,
$$
\end{thm}

\begin{proof}
Consider the problem \eqref{Dual} and the following one, for $\psi\in L^2_\rho(\mathbb R^2)$,
\begin{equation}
\int_{\mathbb R^2}\rho^2u\,\overline\psi\,dx = \mu I\,(\overline{\varphi}(z_1)-\overline{\varphi}(z_2)).
\end{equation}
where $\varphi$ is the solution of Problem \eqref{AdjointL}. Then
\begin{equation}
\int_{\mathbb R^2}\rho^2(u\eps-u)\overline\psi\,dx =
\mu I\Big(\frac 1{|\Omega\eps_1|}\int_{\Omega\eps_1}\varphi\eps\,dx-\varphi(z_1)\Big) -
\mu I\Big(\frac 1{|\Omega\eps_2|}\int_{\Omega\eps_2}\varphi\eps\,dx-\varphi(z_2)\Big).
\label{eq1-1}
\end{equation}
Since $\varphi\in H^2(B)\subset C^{0,\alpha}(\overline B)$ for all $\alpha$ with $0<\alpha<1$
(see \cite{B} for instance) and all compact subsets $B$ of $\mathbb R^2$, we have for $k=1,2$,
\begin{align}
\Big|\frac 1{|\Omega_k\eps|}\int_{\Omega_k\eps}\varphi(x)\,dx - \varphi(z_k)\Big|
&\le \frac 1{|\Omega_k\eps|}\int_{\Omega_k\eps}|\varphi(x) - \varphi(z_k)|\,dx\nonumber\\
&\le C\,\frac 1{|\Omega_k\eps|}\int_{\Omega_k\eps}|x-z_k|^\alpha\,dx\nonumber\\
&\le C\,\epsilon^\alpha.\label{eq2}
\end{align}
Furthermore, we have from \eqref{ErrorPhi}, the imbedding $H^2(B)\subset C^0(\overline B)$
and the mean value theorem,
\begin{align}
\frac 1{|\Omega_k\eps|}\Big|\int_{\Omega_k\eps}(\varphi\eps-\varphi)\,dx\Big|
&\le C_1\,\|\varphi\eps-\varphi\|_{C^0(B)}\nonumber\\
&\le C_2\,\|\varphi\eps-\varphi\|_{H^2(B)}\nonumber\\
&\le C_3\epsilon\,\|\rho\psi\|_{L^2(\mathbb R^2)}.\label{eq3}
\end{align}
Recalling \eqref{eq1-1} and using \eqref{eq2}, \eqref{eq3}, we get
$$
\lim_{\epsilon\to 0}\int_{\mathbb R^2}(u\eps-u)\,\rho^2\overline\psi\,dx = 0\qquad\forall\ \psi\in L^2_\rho(\mathbb R^2).
$$
The sequence $(u\eps)$ converges then weakly to $u$ in $L^2_\rho(\mathbb R^2)$.
To obtain the strong convergence of $u\eps$, we choose $\psi=(u\eps-u)\in L^2_\rho(\mathbb R^2)$ in \eqref{eq1-1}.
We have by using again \eqref{eq2}, \eqref{eq3},
\begin{align*}
\|\rho(u\eps-u)\|^2_{L^2(\mathbb R^2)} &\le \mu I \sum_{k=1}^2 \Big|\int_{\Omega_k\eps}(\varphi\eps-\varphi)\,dx\Big|
+ \mu I \sum_{k=1}^2\left|\frac 1{|\Omega_k\eps|}\int_{\Omega_k\eps}\varphi\,dx-\varphi(z_k)\right|\\
&\le C_4\,\epsilon + C_5\,\epsilon^\alpha \le C\,\epsilon^\alpha.
\end{align*}
\end{proof}

\section{Sharper convergence results}

The convergence result obtained in the previous section can be improved, as we shall show hereafter,
by using the technique of renormalized solutions for elliptic equations following Boccardo -- Gallou\"et \cite{BG}
and Murat \cite{Murat}. To simplify the settings, we shall sometimes resort to writing
Problem \eqref{Pb-Ie} as a system of two coupled equations
involving real valued unknowns. Let us denote, for a complex number $z$, by $z_R$ and $z_I$
its real and imaginary parts respectively. Equation \eqref{Pb-Ie} can be written:
\begin{alignat}{2}
&-\Delta u\eps_R - \beta\sum_{k=1}^2\chi_k\eps\,(u\eps_I-\wt u\eps_{k,I})
- \beta\chi_0\, u\eps_I = \mu I\Big(\frac{\chi_1\eps}{|\Omega_1\eps|}
-\frac{\chi_2\eps}{|\Omega_2\eps|}\Big)&&\quad\text{in }\mathbb R^2,\label{Pb-Ie-R-1}\\
&-\Delta u\eps_I + \beta\sum_{k=1}^2\chi_k\eps\,(u\eps_R-\wt u\eps_{k,R})
+ \beta\chi_0\, u\eps_R = 0&&\quad\text{in }\mathbb R^2,\label{Pb-Ie-R-2}\\
&u_R\eps(x) = \alpha_R + O(|x|^{-1})&&\quad\text{as } |x|\to\infty,\label{Pb-Ie-R-3}\\
&u_I\eps(x) = \alpha_I + O(|x|^{-1})&&\quad\text{as } |x|\to\infty.\label{Pb-Ie-R-4}
\end{alignat}

We start by deriving $L^2$ and $L^1$ uniform estimates.

\begin{lem}
\label{ErrEst-L2-A}
We have the estimates:
\begin{align}
&\|\rho\,u\eps\|_{L^2(\mathbb R^2)} + \epsilon^{-\frac 12}\sum_{k=1}^2\|u\eps-\wt u_k\eps\|_{L^2(\Omega_k\eps)}\le C,\label{Estim-L2}\\
&\|\rho^2 u\eps\|_{L^1(\mathbb R^2)} + \epsilon^{-\frac 32}\sum_{k=1}^2\|u\eps-\wt u_k\eps\|_{L^1(\Omega_k\eps)}\le C.\label{Estim-L1}
\end{align}
\end{lem}

\begin{proof}
The estimate on $\|\rho\,u\eps\|_{L^2(\mathbb R^2)}$ is obtained from Theorem \ref{ErrEst-L2} and from the
fact that $\rho\, u\in L^2(\mathbb R^2)$.
Next, The H\"older's inequality gives
$$
\int_{\mathbb R^2}\rho^2|u\eps|\,dx \le \bigg(\int_{\mathbb R^2}|\rho\, u\eps|^2\,dx\bigg)^{\frac 12}\,
\bigg(\int_{\mathbb R^2}\rho^2\,dx\bigg)^{\frac 12} \le C_2\bigg(\int_{\mathbb R^2}|\rho\, u\eps|^2\,dx\bigg)^{\frac 12}.
$$
Using a variational formulation of Problem \eqref{Pb-Ie}, we then obtain the bound
$$
\int_{\mathbb R^2}|\nabla u\eps|^2\,dx \le C_1\sum_{k=1}^2\frac 1{|\Omega_k\eps|}\|u\eps\|_{L^1(\Omega_k\eps)}
\le C_2\,\epsilon^{-1}.
$$
The Poincar\'e-Wirtinger inequality yields for $k=1,2$,
$$
\int_{\Omega_k\eps}|u\eps-\wt u_k\eps|^2\,dx \le C_1\,\epsilon^2\int_{\mathbb R^2}|\nabla u\eps|^2\,dx \le C_2\,\epsilon.
$$
Again, the Cauchy-Schwarz inequality gives the $L^1$--estimate:
$$
\int_{\Omega_k\eps}|u\eps-\wt u_k\eps|\,dx \le |\Omega_k\eps|^{\frac 12}\Big(\int_{\Omega_k\eps}|u\eps-\wt u_k\eps|^2\,dx\Big)^{\frac 12}
\le C\,\epsilon^{\frac 32}.
$$
\end{proof}


%

We now need a technical result before proving a convergence result. The result, which is a variant of
the Poincar\'e--Wirtinger inequality, can be established by an analogous proof.

\begin{lem}
\label{PW-p}
There exists a constant $C$ such that
$$
\|v\|_{L^p(B)} \le C\,\|\nabla v\|_{L^p(B)^2}\qquad\forall\ v\in W^{1,p}(B) \text{ with }\int_{\Omega_0} v\,dx=0,
$$
where $1\le p<\infty$ and $B$ is any compact subset of $\mathbb R^2$ that contains $\Omega\eps$.
\end{lem}

\begin{thm}
The sequence $(u\eps)$ converges weakly in $W^{1,p}(B)$, $1\le p<2$, toward the unique solution $u$
of Problem \eqref{Pb-Lim} in each ball $B$ containing $\overline\Omega\eps$.
\end{thm}

\begin{proof}
For an integer $m$, we define a subset $B_m\eps$ of $B$ by
$$
B\eps_m = \{x\in B;\ 2^m\le \max\{|u_R\eps(x)|,|u_I\eps(x)|\}\le 2^{m+1}\}.
$$
Let $\psi_m$ stand for the truncature function defined by
$$
\psi_m(s) := \begin{cases}
0 &\text{if}\quad 0\le s\le 2^m,\\
s-2^m &\text{if}\quad 2^m\le s\le 2^{m+1},\\
2^m &\text{if}\quad 2^{m+1}\le s,
\end{cases}
$$
extended to $\mathbb R$ by oddity. Multiplying Equation \eqref{Pb-Ie-R-1} by $\psi_m(u_R\eps)$ and
Equation \eqref{Pb-Ie-R-2} by $\psi_m(u_I\eps)$,
integrating on $\mathbb R^2$, using the Green formula and summing up, we get
\begin{multline}
\int_{\mathbb R^2}\psi'_m(u_R\eps)\,|\nabla u_R\eps|^2\,dx +
\int_{\mathbb R^2}\psi'_m(u_I\eps)\,|\nabla u_I\eps|^2\,dx -
\beta\int_{\mathbb R^2}\chi_0 u_I\eps\psi_m(u_R\eps)\,dx \\
+\beta\int_{\mathbb R^2}\chi_0 u_R\eps\psi_m(u_I\eps)\,dx
- \beta\sum_{k=1}^2\int_{\mathbb R^2}\chi_k\eps(u_I\eps-\wt u_{I,k}\eps)\psi_m(u_R\eps)\,dx \\
+ \beta\sum_{k=1}^2\int_{\mathbb R^2}\chi_k\eps(u_R\eps-\wt u_{R,k}\eps)\psi_m(u_I\eps)\,dx
= \int_{\mathbb R^2}r\eps\psi_m(u_R\eps)\,dx,\label{ident-1}
\end{multline}
where
$$
r\eps = \mu I\Big(\frac{\chi_1\eps}{|\Omega_1\eps|}-\frac{\chi_2\eps}{|\Omega_2\eps|}\Big).
$$
Note that we have
\begin{equation}
\|r\eps\|_{L^1(\mathbb R^2)}\le C.\label{bound-r}
\end{equation}
Since $\psi'_m\ge 0$ and $|\psi_m(u\eps)|\le 2^m$, we have by using \eqref{Estim-L1},
\begin{equation}
\frac 1{2^m}\int_{B_m\eps}|\nabla u\eps|^2\,dx \le C,\label{Est-1}
\end{equation}
where $C$ is independent of $\epsilon$ and $m$. Let $p$ denote a real number with
$1<p<2$. We have from the H\"older inequality
\begin{equation}
\int_{B_m\eps}|\nabla u\eps|^p\,dx \le \Big(\int_{B_m\eps}|\nabla u\eps|^2\,dx\Big)^{\frac p2}
|B\eps_m|^{1-\frac p2}.\label{Est-1a}
\end{equation}
Since $|u\eps|\ge 2^m$ on $B_m\eps$, we have by using the H\"older inequality,
$$
|B\eps_m|\le \frac 1{2^m}\int_{B\eps_m}|u\eps|\,dx
\le \frac 1{2^m}\Big(\int_{B\eps_m}|u\eps|^s\,dx\Big)^{\frac 1s}|B_m\eps|^{\frac 1{s'}}
$$
for all $s,s'\ge 1$ with $1/s + 1/s'=1$. Hence
$$
|B_m\eps|\le \frac 1{2^{ms}}\Big(\int_{B\eps_m}|u\eps|^s\,dx\Big).
$$
Using \eqref{Est-1a} and \eqref{Est-1} yields then
$$
\int_{B\eps_m}|\nabla u\eps|^p\,dx \le \frac C{2^{m(s(1-p/2)-p/2)}}
\Big(\int_{B\eps_m}|u\eps|^s\,dx\Big)^{1-\frac p2}.
$$
We choose here $s>p/(2-p)$ so that $s(1-p/2)-p/2>0$. Therefore
\begin{equation}
\sum_{m\ge 0}\int_{B_m\eps}|\nabla u\eps|^p\,dx \le
C\sum_{m\ge 0}\frac 1{2^{m(s(1-p/2)-p/2)}}
\Big(\int_{B_m\eps}|u\eps|^s\,dx\Big)^{1-\frac p2}.\label{Est-2}
\end{equation}
From the discrete H\"older inequality
$$
\sum_m a_mb_m \le \Big(\sum_m a_m^r\Big)^{\frac 1r}\Big(\sum_m b_m^{r'}\Big)^{\frac 1{r'}}
\quad\text{for }r,r'\ge 1,\ \frac 1r+\frac 1{r'}=1,
$$
Inequality \eqref{Est-2} yields
$$
\sum_{m\ge 0}\int_{B_m\eps}|\nabla u\eps|^p\,dx
\le C\,\Big(\sum_{m\ge 0}\frac 1{2^{mr(s(1-p/2)-p/2)}}\Big)^{\frac 1r}
\Big(\sum_{m\ge 0}\Big(\int_{B_m\eps}|u\eps|^s\,dx\Big)^{r'(1-p/2)}\Big)^{\frac 1{r'}}.
$$
Choosing $r'=2/(2-p)$, we obtain
\begin{equation}
\sum_{m\ge 0}\int_{B\eps_m}|\nabla u\eps|^p\,dx \le C\Big(\sum_{m\ge 0}\int_{B_m\eps}|u\eps|^s\,dx\Big)^{1-\frac p2}.
\label{Est-3}
\end{equation}
We next define
$$
\wt B\eps = \{x\in B;\ 0\le \max\,\{|u_R\eps(x)|,|u_I\eps(x)|\}\le 1\},
$$
which clearly implies
$$
B = \wt B\eps\cup (\bigcup_{m\ge 0}B_m\eps).
$$
In order to estimate $u\eps$ in $W^{1,p}(\wt B\eps)$, we define the truncation function
$$
T(s) = \begin{cases}
\phantom{-}1 &\text{if }s\ge 1\\
\phantom{-}s &\text{if }-1\le s\le 1\\
-1 &\text{if }s\le -1.
\end{cases}
$$
Multiplying Equation \eqref{Pb-Ie-R-1} by $T(u_R\eps)$, Equation \eqref{Pb-Ie-R-2} by $T(u_I\eps)$,
integrating on $\mathbb R^2$, using the Green formula and summing up, we obtain
\begin{multline*}
\int_{\mathbb R^2}T'(u_R\eps)|\nabla u_R\eps|^2\,dx +
\int_{\mathbb R^2}T'(u_I\eps)|\nabla u_I\eps|^2\,dx -
\beta\int_{\mathbb R^2}\chi_0 u_I\eps\, T(u_R\eps)\,dx\\
+\beta\int_{\mathbb R^2}\chi_0 u_R\eps\, T(u_I\eps)\,dx
- \beta\sum_{k=1}^2\int_{\mathbb R^2}\chi_k\eps(u_I\eps-\wt u_{I,k}\eps)\, T(u_R\eps)\,dx \\
+ \beta\sum_{k=1}^2\int_{\mathbb R^2}\chi_k\eps(u_R\eps-\wt u_{R,k}\eps)\, T(u_I\eps)\,dx
= \int_{\mathbb R^2}r\eps\, T(u_R\eps)\,dx.
\end{multline*}
Using \eqref{Estim-L1}, the bound \eqref{bound-r} and the properties $|T(s)|\le 1$, $T'\ge 0$, we deduce
\begin{align*}
\int_{\wt B\eps}|\nabla u_R\eps|^2\,dx + \int_{\wt B\eps}|\nabla u_I\eps|^2\,dx
&\le \|r\eps\|_{L^1(\mathbb R^2)} + \beta\,\Big(\|u_I\eps\|_{L^1(\Omega_0)} + \|u_R\eps\|_{L^1(\Omega_0)}\\
&\qquad + \sum_{k=1}^2\|u_R\eps-\wt u\eps_{R,k}\|_{L^1(\Omega_k\eps)}
+ \sum_{k=1}^2\|u_I\eps-\wt u\eps_{I,k}\|_{L^1(\Omega_k\eps)}\Big)\\
&\le C.
\end{align*}
This yields
\begin{equation}
\int_{\wt B\eps}|\nabla u\eps|^p\,dx \le C.\label{Est-4}
\end{equation}
Combining \eqref{Est-4} and \eqref{Est-3}, we have then in particular
\begin{equation}
\int_B|\nabla u\eps|^p\,dx \le C\,\left(1+\Big(\int_B|u\eps|^s\,dx\Big)^{1-\frac p2}\right)
\quad\text{for }s>\frac p{2-p}.\label{est-p-1}
\end{equation}
We use successively the Gagliardo-Nirenberg (see Friedman \cite{Friedman}, p.~27) and Lemma \ref{PW-p}
to get
$$
\Big(\int_B |u\eps|^s\,dx\Big)^{\frac 1s} \le C\Big(\int_B|\nabla u\eps|^p\,dx\Big)^{\frac\lambda p}
\Big(\int_B|u\eps|\,dx\Big)^{1-\lambda},
$$
with $0\le\lambda\le 1$ and such that
$$
\lambda = \frac{1- \frac 1s}{\frac 32-\frac 1p}.
$$
Using \eqref{Estim-L1} yields
$$
\int_B |u\eps|^s\,dx\le C\,\Big(\int_B|\nabla u\eps|^p\Big)^{\frac{\lambda s}p},
$$
where $C=C(B)$.
Whence, from \eqref{est-p-1},
\begin{equation}
\int_B|\nabla u\eps|^p\,dx \le C\,\left(1+\Big(\int_B|\nabla u\eps|^p\,dx\Big)^{\frac{\lambda s(1-p/2)}{p}}\right),
\label{est-p-2}
\end{equation}
for all $s>p/(2-p)$ and $0\le\lambda\le 1$. Let us choose for $s$ the value
$(1+p)/(2-p)$ that yields
$$
\frac{\lambda s}p\Big(1-\frac p2\Big) = \frac{2-p}{3p-2} < 1 \quad\text{ for }1<p<2.
$$
We then deduce from \eqref{est-p-2} the bound
$$
\int_B|\nabla u\eps|^p\,dx \le C,
$$
with $C=C(B)$.
Therefore, the sequence $(u\eps)$ is bounded in $W^{1,p}(B)$ for all balls $B$ that contain $\overline\Omega\eps$.
From this, we deduce that a subsequence of $(u\eps)$, still denoted by $(u\eps)$, satisfies
$$
u\eps\rightharpoonup u^*\qquad\text{in }W^{1,p}(B).
$$
From the compactness of the imbedding $W^{1,p}(B)\subset L^q(B)$ for $1\le q<2p/(2-p)$, we have
$$
u\eps\to u^*\qquad\text{in }L^q(B)\quad\text{for }1\le q<\frac{2p}{2-p}.
$$
Theorem \ref{ErrEst-L2} implies $u^*=u$. Thus, the subsequence of $(u^\epsilon)$ converges
strongly to $u$ in $L^q(B)$.
Let us show that the convergence to the solution of \eqref{Pb-Lim} takes place in $W^{1,p}(B)$--weak for all bounded
balls $B$ of $\mathbb R^2$. We have from \eqref{Pb-Ie} for all
$\varphi\in W^{1,p'}(B)$ extended by zero outside $B$, with $1/p+1/p'=1$,
$$
\int_B\nabla u\eps\cdot\nabla\overline\varphi\,dx + i\beta\sum_{k=1}^2\int_{\Omega_k\eps}(u\eps-\wt u\eps_k)\overline\varphi
+ i\beta\int_{\Omega_0}u\eps\overline\varphi\,dx
= \frac{\mu I}{|\Omega_1\eps|}\int_{\Omega_1\eps}\overline\varphi\,dx
- \frac{\mu I}{|\Omega_2\eps|}\int_{\Omega_2\eps}\overline\varphi\,dx.
$$
We have
\begin{equation}
\label{lim-1}
\int_B\nabla u\eps\cdot\nabla\overline\varphi\,dx\to\int_B\nabla u\cdot\nabla\overline\varphi\,dx.
\end{equation}
Next, using \eqref{Estim-L2}, we have for $k=1,2$,
$$
\big|\int_{\Omega_k\eps}(u\eps-\wt u\eps_k)\overline\varphi\,dx\Big|\le
\|u\eps-\wt u\eps_k\|_{L^2(\Omega_k\eps)}\,\|\varphi\|_{L^2(\Omega_k\eps)}\le C\epsilon^{\frac 12}\,\|\varphi\|_{L^2(\Omega_k\eps)}.
$$
Therefore
\begin{equation}
\label{lim-2}
\int_{\Omega_k\eps}(u\eps-\wt u\eps_k)\,\overline\varphi\,dx\to 0.
\end{equation}
For the term involving $\Omega_0$, we deduce from Lemma \ref{ErrEst-L2-A},
\begin{equation}
\label{lim-3}
\int_{\Omega_0}u\eps\,\overline\varphi\,dx\to \int_{\Omega_0} u\overline\varphi\,dx.
\end{equation}
Finally, since $p'>2$, then we have the imbedding of $W^{1,p'}(B)$ into $C^0(\overline B)$,
which implies
\begin{equation}
\label{lim-4}
\frac{\mu I}{|\Omega_1\eps|}\int_{\Omega_1\eps}\overline\varphi\,dx
- \frac{\mu I}{|\Omega_2\eps|}\int_{\Omega_2\eps}\overline\varphi\,dx
\to \mu I\overline\varphi(z_1) - \mu I\overline\varphi(z_2).
\end{equation}
Collecting \eqref{lim-1}--\eqref{lim-4}, we find for $u$ the equation
$$
\int_B\nabla u\cdot\nabla\overline\varphi\,dx + i\beta\int_{\Omega_0}u\overline\varphi\,dx
=\mu I\,(\overline\varphi(z_1) - \overline\varphi(z_2)).
$$
This implies that $u$ satisfies the first equation of Problem \eqref{Pb-Lim} on $B$.
Thanks to Lemma \ref{Uniqueness}, the whole sequence $(u\eps)$ converges to $u$ weakly in $W^{1,p}(B)$ and strongly
in $L^q(B)$, for $1\le q\le 2p/(2-p)$.
\end{proof}

Let us conclude by some remarks:

\begin{enumerate}
\item It is clear that the analysis carried out in this paper can be easily extended
to the case where the physical properties $\mu$ and $\sigma$ are not constant. We shall
however assume, in this case, that the magnetic permeability is a $W^{1,\infty}$ function.
This is necessary for $H^2$ regularity results.
\item The obtained results are generalizable to an arbitrary number of (``thick" or ``thin") conductors.
\item In the particular case where no ``thick" conductor is present (\emph{i.e.} $\Omega_0=\emptyset$),
the limit problem becomes
$$
-\Delta u = \mu I(\delta_{z_1}-\delta_{z_2})\qquad\text{in }\mathbb R^2.
$$
Clearly, the solution of this equation is given by
$$
u(x) = \frac{\mu I}{2\pi}\,\log\frac{|x-z_2|}{|x-z_1|},\qquad x\in\mathbb R^2.
$$

\end{enumerate}

\bibliographystyle{unsrt}

%
%
\end{document}